\documentclass[final]{siamltex}

\usepackage[parfill]{parskip}    
\usepackage{amsmath,amssymb,latexsym,enumerate,mathrsfs}
\usepackage{hyperref}
\usepackage{array}
\usepackage{times}
\usepackage[]{mdframed}
\usepackage{epstopdf}
\usepackage{subfig}
\usepackage{graphicx}
\usepackage{hyperref} 
\usepackage{multicol}
\usepackage{framed}
\usepackage{moreverb}
\usepackage{marvosym}
\usepackage{color,soul}
\definecolor{lightblue}{rgb}{.90,.95,1}
\sethlcolor{yellow}

\def\qed{\hfill $\Box$}

\title{A successive approximation method in functional spaces for hierarchical optimal control problems and its application to learning}

\author{Getachew K. Befekadu}

\begin{document}
\maketitle

\renewcommand{\thefootnote}{\arabic{footnote}}

\begin{abstract}
In this paper, we consider a class of learning problem of point estimation for modeling high-dimensional nonlinear functions, whose learning dynamics is guided by model training dataset, while the estimated parameter in due course provides an acceptable prediction accuracy on a different model validation dataset. Here, we establish an evidential connection between such a learning problem and a hierarchical optimal control problem that provides a framework how to account appropriately for both generalization and regularization at the optimization stage. In particular, we consider the following two objectives that the learning process is expected to accomplish: (i) The first one being a {\it controllability-type problem}, i.e., a task related to generalization, which consists of guaranteeing the estimated parameter to reach a certain target set at some fixed final time, where such a target set is associated with model validation dataset. (ii) The second one is a {\it regularization-type problem} that ensures the trajectory of the estimated parameter to satisfy some regularization property over a certain finite time interval. First, we partition the control into two control strategies that are compatible with two abstract agents, namely, a ``{\it leader}'' (which is responsible for the controllability-type problem) and that of a ``{\it follower}'' (which is associated with the regularization-type problem). Then, using the notion of {\it Stackelberg's optimization}, we provide conditions on the existence of admissible optimal control strategies for such a hierarchical optimal control problem under which the ``{\it follower}'' is required to respond optimally to the strategy of the ``{\it leader},'' so as to achieve the overall objectives that ultimately leading to an optimal parameter estimate. Moreover, we provide a nested algorithm, which is arranged in a hierarchical structure-based on successive approximation methods in functional spaces, for numerically solving the corresponding hierarchical optimal control problem. Finally, as part of this work, we also present some numerical results for a typical case of nonlinear regression problem.
\end{abstract}

\section{Introduction} \label{S1}
In this paper, we consider a class of learning problem of point estimation for modeling high-dimensional nonlinear functions, whose learning dynamics is guided by model training dataset, while the estimated parameter in due course provides an acceptable prediction accuracy on a different model validation dataset. Here, we establish an evidential connection between such a learning problem and a hierarchical optimal control problem that provides a framework how to account appropriately for both generalization and regularization at the optimization stage. In particular, we consider the following two objectives that the learning process is expected to accomplish: (i) The first one being a {\it controllability-type objective}, i.e., a task related to generalization, which consists of guaranteeing the estimated parameter to reach a certain target set at some fixed final time, where such a target set is associated with model validation dataset. (ii) The second one is a {\it regularization-type objective} that ensures the trajectory of the estimated parameter to satisfy some regularization property over a certain finite time interval. First, we partition the control into two control strategies that are compatible with two abstract agents, namely, a ``{\it leader}'' (which is responsible for the {\it controllability-type problem}) and that of a ``{\it follower}'' (which is associated with the {\it regularization-type problem}). Then, using the notion of {\it Stackelberg's optimization} (e.g., see \cite{r1}, \cite{r2}, \cite{r3a} or \cite{r3}) we provide conditions on the existence of admissible optimal control strategies for such a hierarchical optimal control problem under which the ``{\it follower}'' is required to respond optimally to the strategy of the ``{\it leader},'' so as to achieve the overall objectives that ultimately leading to an optimal parameter estimate. Here, it is worth remarking that the proofs for such optimality conditions are a consequence of the {\it Pontryagin's maximum principle} in functional spaces (e.g., see \cite{r4} for additional discussions on the first-order necessary optimality conditions based on the {\it Pontryagin's maximum principle}), for which we also present a nested algorithm, which is arranged in a hierarchical structure-based on successive approximation methods in functional spaces, for solving numerically the corresponding hierarchical optimal control problem. In this paper, as part of this work, we also present some numerical results for a typical case of nonlinear regression problem.

The remainder of this paper is organized as follows. In Section~\ref{S2}, we provide a formal problem statement along with some core concepts that will allow us to establish an evidential connection between a class of learning problem and that of a hierarchical optimal control problem. In Section~\ref{S3}, we present our main results, where we provide the optimality conditions for such an optimal control problem that ultimately leading to an optimal parameter estimate. Here, we also provide two propositions, whose proofs are a consequence of the {\it Pontryagin's maximum principle} in functional spaces. Moreover, we present a nested algorithm, which is arranged in a hierarchical structure-based on successive approximation methods in functional spaces, for numerically solving the corresponding optimal control problem. In Section~\ref{S4}, we present numerical results for a typical nonlinear regression problem, and Section~\ref{S4} contains some concluding remarks. For the sake of readability, all proofs are presented in the Appendix section.

\section{Problem formulation} \label{S2} 
In this section, we provide a formal problem statement, where we establish an evidential connection between a class of learning problem and a hierarchical optimal control problem that provides a mathematical construct how to account appropriately for both generalization and regularization. In particular, the problem statement consists of the following core ideas:
\begin{enumerate} [(i).]
\item {\it Datasets}: We are given two datasets, i.e., $\mathcal{Z}^{(k)} = \bigl\{ (x_i^{(k)}, y_i^{(k)})\bigr\}_{i=1}^{m_k}$, each with data size of $m_k$, for $k=1, 2$. The datasets $\bigl\{ \mathcal{Z}^{(1)} \bigr\}_{k=1}^{2}$ may be generated from a given original dataset $\mathcal{Z}^{(0)} =\bigl\{ (x_i^{(0)}, y_i^{(0)})\bigr\}_{i=1}^{m_0}$ by means of bootstrapping with/without replacement. Here, we assume that the first dataset $\mathcal{Z}^{(1)}=\bigl\{ (x_i^{(1)}, y_i^{(1)})\bigr\}_{i=1}^{m_1}$ will be used for model training purpose, while the second dataset $\mathcal{Z}^{(2)} = \bigl\{ (x_i^{(2)}, y_i^{(2)})\bigr\}_{i=1}^{m_2}$ will be used for evaluating the quality of the estimated model parameter.
\item {\it Learning via gradient systems}: We are tasked to find for a parameter $\theta \in \Theta$, from a finite-dimensional parameter space $\mathbb{R}^p$ (i.e., $\Theta \subset \mathbb{R}^p$), such that the function $h_{\theta}(x) \in \mathcal{H}$, i.e., from a given class of hypothesis function space $\mathcal{H}$, describes best the corresponding model training dataset as well as predicts well with reasonable expectation on a different model validation dataset. Here, the search for an optimal parameter $\theta^{\ast} \in \Gamma \subset \mathbb{R}^p$ can be associated with the {\it steady-state solution} to the following gradient system, whose {\it time-evolution} is guided by the model training dataset $\mathcal{Z}^{(1)}$, i.e.,
\begin{align}
 \dot{\theta}(t) = - \nabla J_0(\theta(t),\mathcal{Z}^{(1)}), \quad \theta(0) = \theta_0, \label{Eq2.1}
\end{align}
with $J_0(\theta, \mathcal{Z}^{(1)}) = \frac{1}{m_1} \sum\nolimits_{i=1}^{m_1} {\ell} \bigl(h_{\theta}(x_i^{(1)}), y_i^{(1)} \bigr)$, where $\ell$ is a suitable loss function that quantifies the lack-of-fit between the model and the datasets.
\item {\it Optimal control problem}: Instead of using directly the gradient system in Equation~\eqref{Eq2.1}, we consider the following controlled-gradient system
\begin{align}
 \dot{\theta}^{u}(t) = - \nabla J_0(\theta^{u}(t),\mathcal{Z}^{(1)}) + u(t), \quad \theta^{u}(0) = \theta_0, \label{Eq2.2}
 \end{align}
where $u(t)$ is a control strategy to be determined from a class of admissible control space $\mathcal{U}$, i.e., $\mathcal{U} = L^{\infty}([0,T], U)$ is the space of measurable and essentially bounded functions from $[0,T]$ to $U \subset \mathbb{R}^p$.\footnote{A function $u \colon [0,T] \to U \subset \mathbb{R}^p$ is said to be essentially bounded, if there exists a set $\omega \subset [0,T]$ of measure zero, such that $u$ is bounded on $[0,T] \setminus \omega.$} 
 Then, we search for an optimal control strategy $u^{\ast}(t) \in \mathcal{U}$, so that we would like to achieve the following two objectives:
 \begin{enumerate} [(1).]
\item {\it Controllability-type objective}: Suppose that we are given a target set $\Gamma \in \Theta$ that may depend on the model validation dataset $\mathcal{Z}^{(2)}$, where such a target set can be used as a criterion for evaluating the quality of the estimated model parameter. Then, the final estimated parameter, at some fixed time $T$, is expected to reach this target set $\Gamma$, i.e.,
\begin{align}
 \theta^{u^{\ast}}(T) \in \Gamma, \label{Eq2.3}
 \end{align}
starting from an initial point $\theta^{u}(0) = \theta_0$.
\item {\it Regularization-type objective }: Under certain conditions (see below the general assumptions), we would like to ensure the estimated parameter trajectory $\theta^{u}(t) \in \Theta^{M}$, with $u(t) \in \mathcal{U}$, where $\Theta^{M}$ is a $p$-dimensional manifold, to satisfy some regularization property over a finite time interval $[0, T]$. 
\end{enumerate}
\item {\it General assumptions:} Throughout the paper, we assume the following conditions\footnote{Note that these conditions are sufficient for the existence of a nonempty compact reachable set $\mathcal{R}(\theta_0)$, for some admissible controls on $[0,T]$ that belongs to $\mathcal{U}$, starting from an initial point $\theta^{u}(0)=\theta_0$ (e.g., see \cite{r5} for related discussions on the Filippov's theorem providing a sufficient condition for compactness of the reachable set). Moreover, the controllability-type problem is solvable, if there exists $u^{\ast}(t) \in \mathcal{U}$ such that $\mathcal{R}(\theta^{u^{\ast}}(T)) \cap \Gamma \neq \emptyset$.}
 \begin{enumerate} [(a).]
 \item The set $U \subset \mathbb{R}^{p}$ is compact and the final time $T$ is fixed.
 \item The set $\tilde{\Theta} = \left \{- \nabla J_0(\theta,\mathcal{Z}^{(1)}) + u \,\bigl\vert\, u \in U \right\}$ is convex for every $\theta \in \Theta^{M}$.
 \item For every $u(t) \in \mathcal{U}$, the solution of the controlled-gradient system $\dot{\theta}^{u}(t) = - \nabla J_0(\theta^{u}(t),\mathcal{Z}^{(1)}) + u(t)$, with $\theta^{u}(0) = \theta_0$, is defined on the whole interval $[0,T]$ and belongs to $\Theta^{M}$.
 \end{enumerate}
\end{enumerate}
In order to make the above optimal control problem mathematically more precise, we consider the following hierarchical cost functionals:
\begin{align}
J_1[u] &= \int_0^T \frac{1}{2} \bigl \Vert \theta^{u}(t) \bigr \Vert^2 dt \notag\\
             & \text{s.t.} \quad \theta^{u}(T) \in \Gamma \label{Eq2.4}
\end{align}
and
\begin{align}
J_2[u] = \int_0^T \left \{ \frac{\alpha}{2} \bigl\Vert \theta^{u}(t) \bigr \Vert^2 + \frac{\beta}{2} \bigl\Vert u(t) \bigr\Vert^2 \right\} dt, \quad \alpha > 0, \quad \beta > 0. \label{Eq2.5}
\end{align}
Note that, in general, finding an optimal control strategy $u^{\ast}(t) \in \mathcal{U}$, i.e., a {\it Pareto-optimal solution}, that simultaneously minimizes the above two cost functionals is not an easy problem. However, in the following section, we introduce the notion of {\it Stackelberg's optimal control} (see \cite{r2}, \cite{r3a}, \cite{r3} or \cite{r1} for additional discussions; see also \cite{r3b} in the context of {\it Stackelberg-Nash controllability-type problem}), where we specifically partition the admissible control strategy space $\mathcal{U}$ into two open subspaces $\mathcal{U}_1$ and $\mathcal{U}_2$, with smooth boundaries ({\it up to a set of measurable} $\mathcal{U}$, with $\mathcal{U}_1 \cap \mathcal{U}_2 = \varnothing$), that are compatible with two abstract agents, namely., a ``{\it leader}'' (which is responsible for the controllability-type problem) and that of a ``{\it follower}'' (which is associated with the regularization-type problem).

Then, the above optimal control problem will become solving for a pair of admissible optimal control strategies $(u_1^{\ast}(t),u_2^{\ast}(t)) \in \mathcal{U}_1 \times \mathcal{U}_2$ that steers the controlled-gradient system
\begin{align}
 \dot{\theta}^{u}(t) = - \nabla J_0(\theta^{u}(t),\mathcal{Z}^{(1)}) + u_1(t) \chi_{\mathcal{U}_1} + u_2(t) \chi_{\mathcal{U}_2}, \quad \theta^{u}(0) = \theta_0 \label{Eq2.6}
\end{align}
so as to achieve the overall objectives of Equations~\eqref{Eq2.4} and \eqref{Eq2.5} that ultimately leading to an optimal parameter estimate $\theta^{u^{\ast}}(T) \in \Gamma$. Here, we use the notation $\chi_{\mathcal{U}_i}$, for $i=1,2$, to denote the characteristic function for $\mathcal{U}_i$, where the admissible control $u_i(t)$ is the restriction of the admissible control $u$ to $\mathcal{U}_i$.

\section{Main results} \label{S3}
In this section, we present our main results, where we provide the optimality conditions for the hierarchical optimal control problem under which the ``{\it follower}'' is required to respond optimally to the strategy of the ``{\it leader},'' so as to achieve the overall objectives that ultimately leading to an optimal parameter estimate. In particular, we provide two propositions, whose proofs are application of the {\it Pontryagin's maximum principle} in functional spaces and, for the sake of readability, their proofs are presented in the Appendix section.

\subsection{Optimality statement for the follower} \label{S3.1}
Let $\tilde{u}_1(t) \in \mathcal{U}_1$ be an admissible control strategy for the ``{\it leader}'' and we further assume that it is given. Moreover, if $u_2^{\ast}(t) \in  \mathcal{U}_2$, i.e., the strategy for the ``{\it follower},'' is an optimal solution to Equations~\eqref{Eq2.5}. Then, such an optimal solution is characterized by the following optimality statement.

\begin{proposition}  \label{P1}
Suppose that an admissible control strategy $\tilde{u}_1(t) \in \mathcal{U}_1$ for the ``{\it leader}'' is given. Then, the solution for the optimal control problem w.r.t. the ``{\it follower,}'' i.e.,
\begin{align}
J_2[u_2] &= \int_0^T \left \{ \frac{\alpha}{2} \bigl \Vert \theta^{u_2}(t) \bigr \Vert^2  + \frac{\beta}{2} \bigl \Vert u_2(t) \chi_{\mathcal{U}_2} \bigr \Vert^2 \right\} dt \quad \to \quad \min_{u_2() \in \mathcal{U}_2} \label{Eq3.1}\\
             & \text{s.t.} \notag \\
             \quad \dot{\theta}^{u_2}(t) &= -\nabla J_0(\theta^{u_2}(t), \mathcal{Z}^{(1)}) + \tilde{u}_1(t) \chi_{\mathcal{U}_1} + u_2(t) \chi_{\mathcal{U}_2}, \quad \theta^{u_2}(0) = \theta_0, \label{Eq3.2}
\end{align}
satisfies the following Euler-Lagrange equations:
\begin{enumerate} [(i).]
\item the forward-equation, i.e., the state estimation dynamics\footnote{Notice that the solution for the state estimation trajectory $\theta^{u_2}(\cdot) \in \Theta^{M}$ and it depends on both the strategies of the ``{\it follower}'' and the ``{\it leader}.''}
\begin{align}
 \dot{\theta}^{u_2}(t) &= \frac{\partial H_2(\theta^{u_2}, p_2, u_2)}{\partial p_2} \notag \\
                                  &=-\nabla J_0(\theta^{u_2}(t), \mathcal{Z}^{(1)}) + \tilde{u}_1(t) \chi_{\mathcal{U}_1} + u_2(t) \chi_{\mathcal{U}_2}, \quad \theta^{u_2}(0) = \theta_0, \label{Eq3.3}
\end{align}
\item the backward-equation, i.e., the adjoint state equation w.r.t. the ``{\it follower,}'' 
\begin{align}
  \dot{p}_2(t) = - \frac{\partial H_2(\theta^{u_2}, p_2, u_2)}{\partial \theta^{u_2}}, \quad p_2(T) = 0, \label{Eq3.4}
\end{align}
\item the extremum condition
\begin{align}
  \frac{\partial H_2(\theta^{u_2}, p_2, u_2)}{\partial u_2} = 0 \quad {\text on} \quad [0,T], \label{Eq3.5}
\end{align}
\end{enumerate}
where $H_2$ is the Hamiltonian equation w.r.t. the ``{\it follower}'' and it is given by
\begin{align}
H_2(\theta^{u_2}, p_2, u_2) &= \bigl \langle -\nabla J_0(\theta^{u_2}(t), \mathcal{Z}^{(1)}) + \tilde{u}_1(t) \chi_{\mathcal{U}_1} + u_2(t) \chi_{\mathcal{U}_2}, p_2 \bigr\rangle  \notag \\
& \quad \quad + \frac{\alpha}{2} \bigl\Vert \theta^{u_2}(t) \bigr\Vert^2  + \frac{\beta}{2} \bigl\Vert u_2(t) \chi_{\mathcal{U}_2} \bigr \Vert^2. \label{Eq3.6}
\end{align}
\end{proposition}

Note that that Proposition~\ref{P1} necessitates the following argument, where the ``{\it follower}'' is required to respond optimally to the strategy of the ``{\it leader},'' i.e., the optimal strategy $u_2 \in \mathcal{U}_2$ for the ``{\it follower}'' is a functional on $[0, T]$ that depends on $u_1(\cdot) \in \mathcal{U}_1$. Hence, such a correspondence problem between the ``{\it leader}'' and that of the ``{\it follower}'' admits a unique functional mapping $\mathcal{F} [~]$, i.e.,
\begin{align}
u_2(t) = \mathcal{F}[u_1(t) \chi_{\mathcal{U}_1}] \in \mathcal{U}_2 \quad \left (\text{i.e.,} \quad \mathcal{F} \colon \mathcal{U}_1 \mapsto \mathcal{U}_2 \quad \text{on}\quad [0,T] \right). \label{Eq3.7}
\end{align}
 
\subsection{Optimality statement for the leader} \label{S3.2}
In what follow, we provide an optimality condition for the strategy of the the ``{\it leader},'' when the strategy of the ``{\it follower}'' satisfies the optimality condition of Proposition~\ref{P1}. Here, we restrict our discussion when the final estimated parameter, at some fixed time $T$, satisfies a terminal condition $\Phi(\theta^{u_1}(T),\mathcal{Z}^{(2)}) = z$, where $z$ is a known small positive number and  $\partial \Phi(\theta,\mathcal{Z}^{(2)})/\partial \theta$ satisfies the Lipschitz condition in $\theta$.

\begin{proposition}  \label{P2}
Suppose that Proposition~\ref{P1} holds true.\footnote{i.e., there exists $u_2(t) = \mathcal{F}[u_1(t)\chi_{\mathcal{U}_1}] \in \mathcal{U}_2$, when $u_1(t) \in \mathcal{U}_1$.} Then, the solution for the optimal control problem w.r.t. the ``{\it leader,}'' i.e.,
\begin{align}
J_1[u_1] &= \int_0^T \frac{1}{2} \bigl \Vert \theta^{u_1}(t) \bigr \Vert^2 dt \quad \to \quad \min_{u_1() \in \mathcal{U}_1} \label{Eq3.8}\\
             & \text{s.t.} \quad\quad \Phi(\theta^{u_1}(T),\mathcal{Z}^{(2)}) = z, \label{Eq3.9} \\
             \quad \dot{\theta}^{u_1}(t) &= -\nabla J_0(\theta^{u_1}(t), \mathcal{Z}^{(1)}) + u_1(t) \chi_{\mathcal{U}_1} + \mathcal{F}[u_1(t)\chi_{\mathcal{U}_1}] \chi_{\mathcal{U}_2}, \quad \theta^{u_1}(0) = \theta_0, \label{Eq3.10} 
\end{align}
satisfies the following Euler-Lagrange equations:
\begin{enumerate} [(i).]
\item the forward-equation, i.e., the state estimation dynamics
\begin{align}
 \dot{\theta}^{u_1}(t) &= \frac{\partial H_1(\theta^{u_1}, p_1, u_1)}{\partial p_1} \notag \\
                                  &=-\nabla J_0(\theta^{u_1}(t), \mathcal{Z}^{(1)}) + u_1(t) \chi_{\mathcal{U}_1} + \mathcal{F}[u_1(t)\chi_{\mathcal{U}_1}] \chi_{\mathcal{U}_2}, \quad \theta^{u_1}(0) = \theta_0, \label{Eq3.11}
\end{align}
\item the backward-equation, i.e., the adjoint state equation w.r.t. the ``{\it leader},'' 
\begin{align}
 \dot{p}_1(t) = - \frac{\partial H_1(\theta^{u_1}, p_1, u_1)}{\partial \theta^{u_1}}, \quad p_1(T) = -\frac{\partial \Phi(\theta,\mathcal{Z}^{(1)})} {\partial \theta} \biggr \vert_{\theta=\theta^{u_1}(T)}, \label{Eq3.12}
                                  \end{align}
\item the extremum condition
\begin{align}
 \frac{\partial H_1(\theta^{u_1}, p_1, u_1)}{\partial u_1} = 0 \quad {\text on} \quad [0,T],  \label{Eq3.13}
\end{align}
\end{enumerate}
where $H_1$ is the Hamiltonian equation w.r.t. the ``{\it leader}'' and it is given by
\begin{align}
H_1(\theta^{u_1}, p_1, u_1) &= \bigl \langle -\nabla J_0(\theta^{u_1}(t), \mathcal{Z}^{(1)}) + u_1(t) \chi_{\mathcal{U}_1} + \mathcal{F}[u_1(t)\chi_{\mathcal{U}_1} ] \chi_{\mathcal{U}_2}, p_1 \bigr\rangle  \notag \\
& \quad \quad + \frac{1}{2} \bigl \Vert \theta^{u_1}(t) \bigr \Vert^2. \label{Eq3.14}
\end{align}
\end{proposition}

In the following section, we provide a nested algorithm, based on successive approximation methods, that will numerically solve the optimality conditions stated in the above two propositions, i.e., Propositions~\ref{P1} and \ref{P2}.

\subsection{Algorithm -- Successive approximation method} \label{S3.3}
In what follows, we present a nested algorithm, which is arranged in a hierarchical structure based-on a successive approximation method, for solving numerically the hierarchical optimal control problem in Equations~\eqref{Eq2.4} and \eqref{Eq2.5} subject to the controlled-gradient system of Equation~\eqref{Eq2.2}. The idea behind this nested algorithm is to embed the {\it regularization-type problem} at each time step within the {\it controllability-type problem} and, then solve numerically the the corresponding Euler-Lagrange equations until the extremum condition w.r.t. the ``{\it leader}'' meets some convergence criterion. Here, the interpretation is basically such a pair of optimal admissible control strategies from $\mathcal{U}_1 \times \mathcal{U}_2$ on $[0, T]$ will guarantee to achieve the overall objectives that ultimately leading to an optimal parameter estimate.

Then, the nested algorithm for solving numerically the hierarchical optimal control problem can be summarized as follows:
{\rm \footnotesize

{\bf ALGORITHM:} Nested Algorithm (with Hierarchical Structure) Based-on Successive Approximation Methods
\begin{itemize}
\item[{\bf 0.}] {\bf Initialize:} Start with any admissible control strategy $u_1^{(n-1)}(t) \in \mathcal{U}_1$ (i.e., an initial guess control function on $[0,T]$) for the ``{\it leader}.'' 
\item[{\bf 1.}] {\bf The regularization-type problem:} Choose an admissible control strategy $u_2^{(n-1)}(t) \in \mathcal{U}_2$ (i.e., an initial guess control function on $[0,T]$) for the ``{\it follower}.'' 
\item[{\bf 2.}] Then, using the admissible control pairs $\bigl(u_1^{(n-1)}(t), u_2^{(n-1)}(t)\bigr)$, solve the forward and backward-equations w.r.t. the system dynamics of the ``{\it follower},'' i.e., Equations~\eqref{Eq3.3} and \eqref{Eq3.4}.
\item[{\bf 3.}] Compute the correction term $\delta u_2(\cdot)$ on $[0,T]$ w.r.t. the admissible control $u_2$ of the ``{\it follower}'' (cf. Equation~\eqref{Eq3.5}) using 
 \begin{align*}
  \delta u_2(\cdot) = \gamma_2 \frac{\partial H_2}{\partial u_2} \quad {\text on} \quad [0,T],
\end{align*}
where $\gamma_2 \in (0, 1]$.
 \item[{\bf 4.}] Update the new admissible control for the ``{\it follower}, i.e., $u_2^{(n)}(t)$, using
\begin{align*}
 u_2^{(n)} (t) = u_2^{(n-1)} (t) + \delta u_2(t).
\end{align*}
 \item[{\bf 5.}] {\bf The controllability-type problem:} Using the admissible control strategy pairs $\bigl(u_1^{(n-1)}(t), u_2^{(n)}(t)\bigr)$, with the updated control strategy $u_2^{(n)}(t)$, solve the forward and backward-equations w.r.t. the system dynamics of the ``{\it leader,}'' i.e., Equations~\eqref{Eq3.11} and \eqref{Eq3.12}.
 \item[{\bf 6.}] Compute the correction term $\delta u_1(\cdot)$ on $[0,T]$ w.r.t. the admissible control $u_1$ of the ``{\it leader}'' (cf. Equation~\eqref{Eq3.13}) using 
\begin{align*}
\delta u_1(\cdot) = \gamma_1 \frac{\partial H_1}{\partial u_1} \quad {\text on} \quad [0,T], 
\end{align*}
where $\gamma_1 \in (0, 1]$.
 \item[{\bf 7.}] Update the new admissible control for the ``{\it leader},''  i.e., $u_1^{(n)}(t)$, using
 \begin{align*}
 u_1^{(n)}(t) = u_1^{(n-1)}(t) + \delta u_1(t).
 \end{align*}
\item[{\bf 8.}] With the updated admissible control strategy pairs $\bigl(u_1^{(n)}(t), u_2^{(n)}(t)\bigr)$, repeat Steps $2$ through $7$, until convergence, i.e., $\bigl\Vert \partial H_1/ \partial u_1 \bigr\Vert \le \epsilon_{\rm tol}$, for some error tolerance $\epsilon_{\rm tol} > 0$.
 \item[{\bf 10.}] {\bf Output:} Return the optimal estimated parameter value $\theta^{\ast} = \theta^{u^{(n)}}(T)$, with $u^{(n)}(t) = (u_2^{(n)}(t), u_2^{(n)}(t))$.
\end{itemize}}

Here, it is worth remarking that the above nested algorithm is a hierarchical implementation of a {\it gradient-based optimization method} in functional spaces (e.g., see \cite{r6}, \cite{r7} and \cite{r8}) for additional discussions). Moreover, instead of directly searching for a pair of optimal control strategies from $\mathcal{U}_1 \times \mathcal{U}_2$, we can approximate $u_i(t)$, with fixed $N$, using
\begin{align*}
u_i(t) = \sum\nolimits_{k = 1}^N a_{i,k} \phi_k(t), \quad i =1, 2, 
\end{align*}
where the functions $\phi_k(t)$ on $[0, T]$, for $k=1,2,\ldots, N$, are appropriately chosen basis functions in advance, but the coefficients $a_{i,k}$, for $k=1,2,\ldots, N$ and $i=1,2$, are unknown, so that they can be determined (within the above nested algorithm) to approximate reasonably well the optimal solutions corresponding to the hierarchical optimal control problem (e.g., see \cite{r8} for additional discussions how this type approach can applied for minimizing functionals). 

 \section{Numerical simulation results} \label{S4} 
 In this section, we presented some numerical results based on the data taken from the Michaelis-Menten's classical paper (see \cite{r9}) that deals with the rate of enzyme-catalyzed reaction. Here, we assume a parametrized model, i.e., a simple nonlinear regression function $v = h_{\theta}(w)$, that relates the initial velocity $v$ as a function of the sucrose concentration $w$ from seven different experiments (see Table~\ref{Tb1})
\begin{align*}
v = h_{\theta} (w) = \frac{\theta_0 \, w}{\theta_1 + w},
\end{align*}
where $\theta_0$ and $\theta_1$ are the parameters to be estimated using the nested algorithm of Subsection~\ref{S3.3}. That is, a hierarchical implementation of the successive approximation method for solving numerically the corresponding hierarchical optimal control problem, whose problem formulation and implementation are evidentially related the underlying parameter estimation problem. 
\begin{table}[h]
\begin{center}
{\footnotesize
\caption{Results from the experiments.} \label{Tb1}
  \begin{tabular}{c|c|c}\hline \hline
   &$w$  & $v$ \\
    {\rm Experiments} &{\rm Initial Concentration} & {\rm Initial Velocity} \\ 
    & {\rm  of Sucrose} & \\\hline
    1 &  0.3330 & 3.6360  \\
    2 &  0.1670 & 3.6360 \\
    3 & 0.0833 & 3.2360 \\
    4 &  0.0416 & 2.6660 \\
    5 &  0.0208 & 2.1140 \\
    6 &  0.0104 & 1.4660 \\
    7 &  0.0052 & 0.8661\\
    \hline \hline
  \end{tabular}}
  \end{center}
 \end{table}

\begin{figure}[h]
\begin{center}
 \includegraphics[scale=0.135]{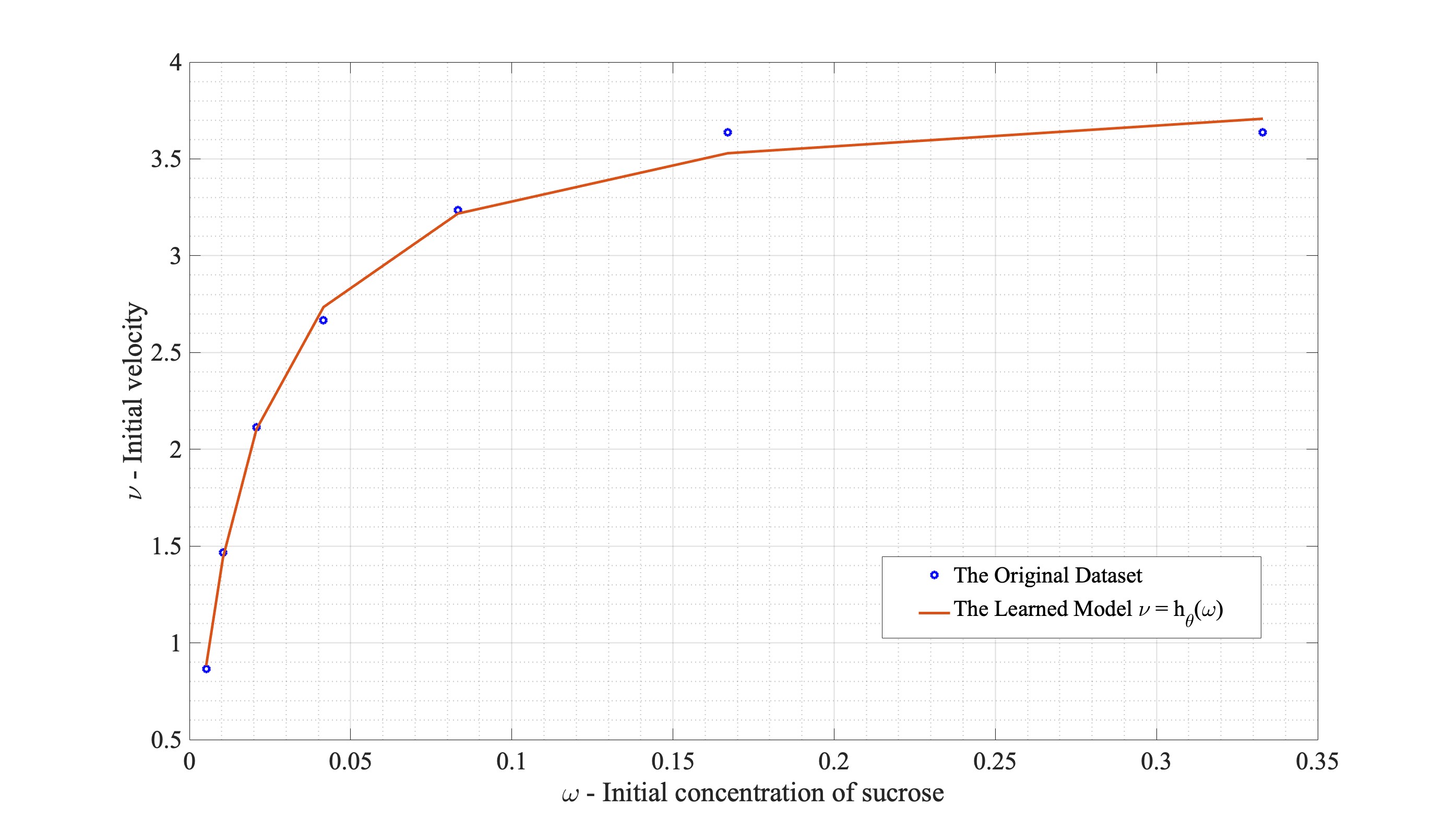}
  \caption{Plots for the original dataset and the learned model $h_{\theta^{\ast}} (w)$.} \label{Fig1}
\end{center}

\end{figure}

\begin{figure}[h]
\begin{center}
   \includegraphics[scale=0.135]{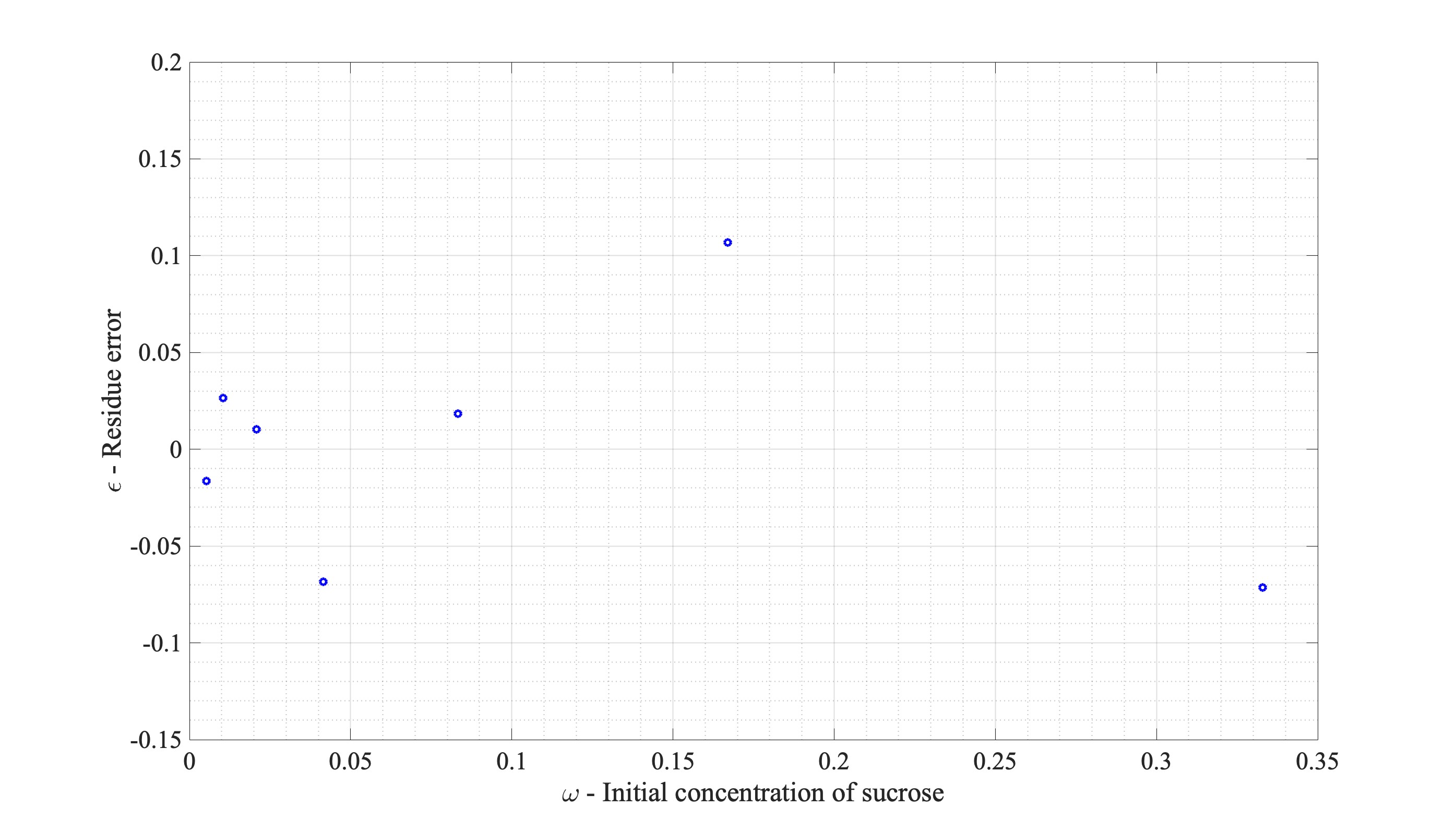}
  \caption{Plot for the residue error $\epsilon$.} \label{Fig2}
\end{center}
\end{figure}

In particular, we used the following problem data specifications:
 \begin{enumerate} [(i)]
 \item The dataset corresponding to the Experiments $\{1,\,3,\, 5, \,7\}$ for the model training purpose, i.e., $\mathcal{Z}^{(1)} = \left\{w_{i_1}^{(1)}, v_{i_1}^{(1)}\right\}_{i_1=1,3,5,7}$, with data size of $m_1=4$; while the dataset corresponding to the Experiments $\{2, \,4,\, 6\}$ for the model validation purpose, i.e., $\mathcal{Z}^{(2)}=\left\{w_{i_2}^{(2)}, v_{i_2}^{(2)}\right\}_{i_2=2,4,6}$, with data size of $m_2=3$.
 \item The gradient system (without a control input term), whose time-evolution is guided by the model training dataset $\mathcal{Z}^{(1)}$ is given by
\begin{align*}
 \dot{\theta}(t) = - \nabla J_0(\theta(t),\mathcal{Z}^{(1)}),
\end{align*}
with $J_0(\theta, \mathcal{Z}^{(1)}) = (1/m_1) \sum\nolimits_{i_1=1}^{m_1} {\ell} \bigl(h_{\theta}(w_{i_1}^{(1)}), v_{i_1}^{(1)} \bigr)$, where $\ell$ is the usual quadratic loss function.
 \item For the controllability-type objective: (a). the parameter values for $\alpha=0.01$ and $\beta=0.1$; and (2) the terminal condition $\Phi(\theta^{u_1}(T),\mathcal{Z}^{(2)}) = z$, with fixed time $T = 1.5$, which is associated with the following loss function
\begin{align*}
\Phi(\theta^{u_1}(T),\mathcal{Z}^{(2)}) =(1/m_2) \sum\nolimits_{i_2=1}^{m_2} {\ell} \bigl(h_{\theta^{u}(T)}(w_{i_2}^{(2)}), v_{i_2}^{(2)} \bigr),
\end{align*}
as a model quality measure, with model accuracy level of $z=0.005$. Note that the regularization-type problem can be embedded within the controllability-type problem at the optimization stage.
 \item An error tolerance $\epsilon_{\rm tol} = 1 \times 10^{-5}$ is used for terminating the nested algorithm; and moreover the algorithm is initialized with constant control values for $u_1$ and $u_2$ on $[0,T]$.
  \end{enumerate}
 Figure~\ref{Fig1} shows the plots for the original dataset and the learned model $h_{\theta} (w)$, with optimal estimated parameter values $\theta_0^{\ast} = 3.9059$ and $\theta_1^{\ast} = 0.0178$.
Moreover, we computed the corresponding sample mean $m=7.9772 \times 10^{-4}$ and standard deviation $s=0.0614$ for the residue error $\varepsilon = v - h(w)$ w.r.t. the original dataset $\mathcal{Z}^{(0)}=\left\{w_{i}, v_{i}\right\}_{i=1}^7$, i.e., the dataset corresponding to all experiments. Moreover, Figure~\ref{Fig2} shows the plot for residue error $\epsilon = \nu - h_{\theta^{\ast}}(\omega)$
\section{Concluding remarks} \label{S5}
In this paper, we considered a class of learning problem and further established an evidential connection with a hierarchical optimal control problem that provides a framework how to account appropriately for both generalization and regularization at the optimization stage. In particular, we considered two objectives: (i). a {\it controllability-type objective}, which is a task related to generalization;  and (ii) a {\it regularization-type objective}, which ensures the trajectory of the estimated parameter to satisfy some regularization property. Moreover, we also provided the optimality conditions for such an optimal control problem, whose proofs are a consequence of the {\it Pontryagin's maximum principle} in functional spaces. Moreover, we presented a nested algorithm, which is arranged in a hierarchical structure-based on successive approximation methods, for solving numerically the corresponding optimal control problem. Finally, as part of this work, we presented some numerical results for a typical nonlinear regression problem.

\section*{Appendix: Proofs}
In this section, we give the proofs for the main results.

{\it Proof of Proposition~\ref{P1}.} Recall that, when $\tilde{u}_1(t) \in \mathcal{U}_2$ is a given admissible control on $[0, T]$. Then, the optimal control problem w.r.t. the ``{\it follower}'' is to minimize the cost functional $J_2[u_2]$, i.e.,
\begin{align}
J_2[u_2] &= \int_0^T \left \{ \frac{\alpha}{2} \bigl \Vert \theta^{u_2}(t) \bigr \Vert^2  + \frac{\beta}{2} \bigl \Vert u_2(t)\chi_{\mathcal{U}_2} \bigr \Vert^2 \right\} dt \quad \to \quad \min_{u_2(\cdot) \in \mathcal{U}_2} \label{Eq5.1}
\end{align}
subject to the evolution of the controlled-gradient system in Equation~\eqref{Eq3.3}.

Introduce the following Lagrange function $\Lambda_2$
\begin{align}
\Lambda_2[p_2, u_2] &= \int_0^T \left \{ \frac{\alpha}{2} \bigl\Vert \theta^{u_2}(t) \bigr\Vert^2  + \frac{\beta}{2} \bigl\Vert u_2(t) \chi_{\mathcal{U}_2}\bigr\Vert^2 \right\} dt \notag\\
                                                        & \quad + \int_0^T \left \langle \dot{\theta}^{u_2}(t) + \nabla J_0(\theta^{u_2}(t), \mathcal{Z}^{(1)}) - \tilde{u}_1(t) \chi_{\mathcal{U}_1} - u_2\chi_{\mathcal{U}_2}, p_2(t) \right \rangle dt, \label{Eq5.2}
\end{align}
where the Lagrange multiplier $p_2(\cdot)$ (i.e., an adjoint state variable w.r.t. the ``{\it follower}'') is a function on $[0, T]$.  Moreover, if we integrate the above equation by parts, then we obtain 
\begin{align}
\Lambda_2[p_2, u_2] &= \left \langle  \theta^{u_2}(T), p_2(T) \right \rangle - \left \langle  \theta^{u_2}(0), p_2(0) \right \rangle \notag\\
                                                        &\quad -  \int_0^T \left \{H_2(\theta^{u_2}(t), p_2(t), u_2(t)) + \left \langle \dot{\theta}^{u_2}(t), p_2(t) \right \rangle \right\} dt,  \label{Eq5.3}
\end{align}
where $H_2$ is the Hamiltonian equation and it is given by
\begin{align}
H_2(\theta^{u_2}, p_2, u_2) &= \bigl \langle -\nabla J_0(\theta^{u_2}(t), \mathcal{Z}^{(1)}) + \tilde{u}_1(t) \chi_{\mathcal{U}_1} + u_2(t) \chi_{\mathcal{U}_2}, p_2 \bigr\rangle  \notag \\
& \quad \quad + \frac{\alpha}{2} \bigl\Vert \theta^{u_2}(t) \bigr\Vert^2  + \frac{\beta}{2} \bigl\Vert u_2(t) \chi_{\mathcal{U}_2} \bigr \Vert^2.  \label{Eq5.4}
\end{align}
Next consider the variation $\delta \Lambda_2$ in $\Lambda_2$, due to the variation of $\delta u_2$ in $u_2$. Note that this variation in the control strategy will further induce a variation of trajectory $\delta \theta^{u_2}$, with $\delta \theta^{u_2}(0) = 0$, since the initial condition $\theta^{u_2}(0)$ is fixed. Then, we have 
\begin{align}
\delta \Lambda_2 = \left \langle \theta^{u_2}(T), p_2(T) \right \rangle &- \int_0^T \left \{ \left \langle \frac{\partial H_2(\theta^{u_2}(t), p_2(t), u_2(t))}{\partial \theta^{u_2}} + \dot{p}_2(t), \delta \theta^{u_2}(t) \right \rangle \right\} dt  \notag\\
                                                        &\quad -  \int_0^T \frac{\partial H_2(\theta^{u_2}(t), p_2(t), u_2(t))}{\partial u_2} \delta u_2(t) dt.  \label{Eq5.5}
\end{align}
Hence, the necessary condition for $\Lambda_2$ to have an extremum is $\delta \Lambda_2 = 0$ for any arbitrary variation of $\delta u_2$, which further implies the following sufficient conditions
\begin{enumerate} [(i).]
\item the backward-equation, i.e., the adjoint state variable $p_2(t)$ satisfies
\begin{align}
\dot{p}_2(t) = - \frac{\partial H_2(\theta^{u_2}(t), p_2(t), u_2(t))}{\partial \theta^{u_2}},  \label{Eq5.6}
\end{align}
\item the forward-equation, i.e., the state estimation dynamics, w.r.t. the ``{\it follower,}'' satisfies Equation~\eqref{Eq3.3}, with a final boundary condition of
\begin{align}
 p_2(T) = 0,  \label{Eq5.7}
\end{align}
\item and further requiring the following extremum condition
\begin{align}
 \frac{\partial H_2(\theta^{u_2}(t), p_2(t), u_2(t))}{\partial u_2} = 0 \quad \text{on} \quad [0,T].  \label{Eq5.8}
\end{align}
\end{enumerate} 
This completes the proof of Proposition~\ref{P1}.\qed

{\it Proof of Proposition~\ref{P2}.}
Recall that, if Proposition~\ref{P1} holds true, then there exists a functional mapping $\mathcal{F}[\,~]$ that transforms a given functional $u_1(\cdot)$ on $[0,T]$ from the class of admissible controls $\mathcal{U}_1$ to another admissible control $u_2(\cdot)$ on $[0,T]$ that belongs $\mathcal{U}_2$, i.e.,
$\mathcal{F}[u_1(t)] \in \mathcal{U}_2$, when $u_1(t) \in \mathcal{U}_1$ .

Note that the optimal control problem w.r.t. the ``{\it leader}'' is to minimize the cost functional $J_1[u_1]$, i.e.,
\begin{align}
J_1[u_1] &= \int_0^T \frac{1}{2} \Vert \theta^{u_1}(t) \Vert^2 dt \quad \to \quad \min_{u_1(\cdot) \in \mathcal{U}_1}  \label{Eq5.9}
\end{align}
subject to the evolution of the controlled-gradient system in Equation~\eqref{Eq3.9}, with a terminal condition $\Phi(\theta^{u_1}(T),\mathcal{Z}^{(1)}) = z$. 

Introduce the following Lagrange function $\Lambda_1$
\begin{align}
\Lambda_1[p_1, u_1] &= \langle \Phi(\theta^{u_1}(T),\mathcal{Z}^{(2)}) - z, p_1^0 \rangle + \int_0^T \frac{1}{2} \Vert \theta^{u_1}(t) \Vert^2 dt \notag\\
                                                        & \quad + \int_0^T \left \langle \dot{\theta}^{u_1}(t) + \nabla J_0(\theta^{u_1}(t), \mathcal{Z}^{(1)}) - u_1(t) \chi_{\mathcal{U}_1} - \mathcal{F}[u_1(t)\chi_{\mathcal{U}_1}], p_1(t) \right \rangle dt,  \label{Eq5.10}
\end{align}
where the Lagrange multiplier $p_1(\cdot)$ (i.e., an adjoint state variable w.r.t. the ``{\it leader}'') is a function on $[0, T]$. Moreover, if we integrate the above equation by parts, then we will obtain the following
\begin{align}
\Lambda_1[p_1, u_1] &= \langle \Phi(\theta^{u_1}(T),\mathcal{Z}^{(2)}) - z, p_1^0 \rangle + \left \langle  \theta^{u_1}(T), p_1(T) \right \rangle - \left \langle  \theta^{u_1}(0), p_1(0) \right \rangle \notag\\
                                                        &\quad -  \int_0^T \left \{H_1(\theta^{u_1}(t), p_1(t), u_1(t)) + \left \langle \dot{\theta}^{u_1}(t), p_1(t) \right \rangle \right\} dt,  \label{Eq5.11}
\end{align}
where $H_1$ is the Hamiltonian equation and it is given by
\begin{align}
H_1(\theta^{u_1}, p_1, u_1) &= \bigl \langle -\nabla J_0(\theta^{u_1}(t), \mathcal{Z}^{(1)}) + u_1(t) \chi_{\mathcal{U}_1} + \mathcal{F}[u_1(t)\chi_{\mathcal{U}_1}], p_1 \bigr\rangle + \frac{1}{2} \bigl \Vert \theta^{u_1}(t) \bigr \Vert^2.  \label{Eq5.12}
\end{align}
Next, consider the variation $\delta \Lambda_1$ in $\Lambda_1$, due to the variation of $\delta u_1$ in $u_1$. Note that this variation in the control strategy will induce a variation of trajectory $\delta \theta^{u_1}$, with $\delta \theta^{u_1}(0) = 0$, since the initial condition $\theta^{u_1}(0)$ is fixed. Then, we have 
\begin{align}
\delta \Lambda_1 =& \left \langle \theta^{u_2}(T), p_1(T) + \frac{\partial \Phi(\theta^{u_1},\mathcal{Z}^{(2)})}{\partial \theta^{u_1}} \biggl \vert_{\theta^{u_1} = \theta^{u_1}(T)} \right \rangle \notag\\
& \quad \quad - \int_0^T \left \{ \left \langle \frac{\partial H_1(\theta^{u_1}(t), p_1(t), u_1(t))}{\partial \theta^{u_1}} + \dot{p}_1(t), \delta \theta^{u_1}(t) \right \rangle \right\} dt  \notag\\
                                                        &\quad\quad \quad  -  \int_0^T \frac{\partial H_1(\theta^{u_1}(t), p_1(t), u_1(t))}{\partial u_1} \delta u_1(t) dt.  \label{Eq5.13}
\end{align}
Hence, the necessary condition for $\Lambda_1$ to have an extremum is $\delta \Lambda_1 = 0$ for any arbitrary variation of $\delta u_1$, which further implies the following sufficient conditions
\begin{enumerate} [(i).]
\item the backward-equation, i.e., the adjoint state variable $p_1(t)$ satisfying
\begin{align}
\dot{p}_1(t) = - \frac{\partial H_1(\theta^{u_1}(t), p_1(t), u_1(t))}{\partial \theta^{u_1}},  \label{Eq5.14}
\end{align}
\item the forward-equation, i.e., the state estimation dynamics, w.r.t. the ``{\it leader,}'' satisfies Equation~\eqref{Eq3.9} with a final boundary condition of
\begin{align}
 p_1(T) = -\frac{\partial \Phi(\theta^{u_1},\mathcal{Z}^{(2)})}{\partial \theta^{u_1}}\biggl \vert_{\theta^{u_1} = \theta^{u_1}(T)}  \label{Eq5.15}
\end{align}
\item and further requiring the following extremum condition
\begin{align}
 \frac{\partial H_1(\theta^{u_1}(t), p_1(t), u_1(t))}{\partial u_1} = 0 \quad \text{on} \quad [0,T].  \label{Eq5.16}
\end{align}
\end{enumerate} 
This completes the proof of Proposition~\ref{P2}.\qed


\begin{thebibliography}{99}

\bibitem{r1}
{H. von Stackelberg}. {\em Marktform und Gleichgewicht}. Springer, Berlin, Germany, 1934.

\bibitem{r2}
{J.L. Lions}. Some remarks on Stackelberg’s optimization. {\em Math. Models Methods Appl. Sci.}, 4, 477--487, 1994.

\bibitem{r3a}
{G.K. Befekadu \& E.L. Pasiliao}. On the hierarchical optimal control of a chain of distributed systems. {\em J. Dynamics and Games}, 2(2), 187--199, 2015.

\bibitem{r3}
{G.K. Befekadu}. Remarks on the hierarchical control problems with model uncertainty. {\em arXiv preprint}, arXiv:1510.03804, 2015.

\bibitem{r4}
{L.S. Pontryagin, V. Boltianski, R. Gamkrelidze, \& E. Mitchtchenko}. {\em The mathematical theory of optimal processes}. John Wiley \& Sons, New York, 1962.

\bibitem{r3b}
{G. Leitmann}. On generalized stackelberg strategies. {\em J. Optim. Theor. Appl.}, 26, 637--643, 1978.

\bibitem{r5}
{D. Liberzon}. {\em Calculus of variations and optimal control theory.} Princeton University Press, Princeton, NJ, 2012.

\bibitem{r6}
{R.H. Moore}. Newton's method and variations. In ``Nonlinear Integral Equations.'' P. M. Anselone (Ed.), pp. 65--98, University of Wisconsin Press, Madison, 1964.

\bibitem{r7}
{M.L. Stein}. On methods for obtaining solutions of fixed-endpoint problems in the calculus of variations. {\em J. Res. Nat. Bur. Standards,} 50(5), 277--297, 1953.

\bibitem{r8}
{S.G. Mikhlin}. {\em Variational methods in mathematical physics}. Macmillan, New York, 1964.

\bibitem{r9}
{L. Michaelis \& M.L. Menten}. Die Kinetik der Invertinwirkung. {\em Biochem. Z.,} 49, 333--369, 1913. 

\end{thebibliography}
\end{document}